\documentclass[11pt]{article}
\usepackage{amssymb}
\usepackage{amscd}
\usepackage{amsmath}
\usepackage{amsfonts}
\usepackage{euscript}
\usepackage{xypic}
\xyoption{all}

\def\scr{\EuScript}

\renewcommand{\thesubsection}{\arabic{section}.\arabic{subsection}}
\newcounter{numero}[subsection]
\renewcommand{\thenumero}{(\thesubsection .\arabic{numero})}
\newenvironment{corolario}{\medskip
\refstepcounter{numero}\noindent {\sc  \thenumero\ Corollary.}\
\it}{\vspace{1ex}\par}

\newenvironment{teorema}{\medskip
\refstepcounter{numero}\noindent {\sc  \thenumero\ Theorem.}\
\it}{\vspace{1ex}\par}

\newenvironment{lema}{\medskip
\refstepcounter{numero}\noindent {\sc  \thenumero\ Lemma.}\
\it}{\vspace{1ex}\par}

\newenvironment{definicion}{\medskip
\refstepcounter{numero}\noindent {\sc  \thenumero\ Definition.}\
\it}{\vspace{1ex}\par}

\newenvironment{proposicion}{\medskip
\refstepcounter{numero}\noindent {\sc  \thenumero\ Proposition.}\
\it}{\vspace{1ex}\par}

\newenvironment{nota}{\medskip
\refstepcounter{numero}\noindent {\sc  \thenumero\ Remark.}\
}{\vspace{1ex}\par}

\newenvironment{ejemplo}{\medskip
\refstepcounter{numero}\noindent {\sc  \thenumero\ Example.}\
}{\vspace{1ex}\par}

\newenvironment{demostracion}{
\noindent {\sc  Proof.}\ }{\hfill Q.E.D.\vspace{1ex}\par}

\newcommand{\numero}{\refstepcounter{numero}\noindent {\sc  \thenumero\ }}

\DeclareMathOperator{\Id}{Id} \DeclareMathOperator{\Perv}{Perv}
 \DeclareMathOperator{\adj}{adj}
\DeclareMathOperator{\cone}{cone} \DeclareMathOperator{\Arr}{Arr}
\DeclareMathOperator{\catder}{D} \DeclareMathOperator{\comp}{C}
\DeclareMathOperator{\Mod}{Mod}

\DeclareMathOperator{\komp}{K} \DeclareMathOperator{\cohom}{H}
\newcommand{\cd}[3]{\catder^{#1}_{#2}(#3)}

\newcommand{\suex}[5]{0 \xrightarrow{} #1 \xrightarrow{#2} #3
                      \xrightarrow{#4} #5 \xrightarrow{} 0}
\newcommand{\suext}[3]{\suex{#1}{}{#2}{}{#3}}
\newcommand{\tri}[3]{#1 \xrightarrow{} #2 \xrightarrow{} #3
                     \xrightarrow{+1}}
\newcommand{\trit}[5]{#1 \xrightarrow{#2} #3
                      \xrightarrow{#4} #5
                      \xrightarrow{} #1 [1]}
\newcommand{\trito}[6]{#1 \xrightarrow{#2} #3
                      \xrightarrow{#4} #5
                      \xrightarrow{#6} #1 [1]}

\newcommand{\cA}{{\cal A}}
\newcommand{\cC}{{\cal C}}
\newcommand{\cB}{{\cal B}}
\newcommand{\calD}{{\cal D}}

\newcommand{\calF}{{\cal F}}

\newcommand{\calL}{{\cal L}}

\newcommand{\calO}{{\cal O}}

\newcommand{\C}{{\mathbb C}}
\newcommand{\F}{{\mathbb F}}

\newcommand{\G}{{\mathbb G}}
\newcommand{\Z}{{\mathbb Z}}
\newcommand{\RR}{{\mathbb R}}
\newcommand{\Q}{{\mathbb Q}}
\newcommand{\QQ}{{\mathbb Q}}
\newcommand{\lra}{\longrightarrow}
\newcommand\Hom{\mbox{\rm Hom}}
\DeclareMathOperator{\coker}{coker}
\DeclareMathOperator{\Img}{Img}

\newcommand\OO{{\scr O}}
\newcommand\gA{\mathfrak{A}}
\newcommand\gB{\mathfrak{B}}
\newcommand\gC{\mathfrak{C}}
\newcommand\gT{\mathfrak{T}}

\newcommand\Omk{\Omega_{\komp}}
\newcommand\Omd{\Omega_{\catder}}
\newcommand\os{\overline{s}}
\def\ot{\overline{t}}
\newcommand\ozeta{\overline{\zeta}}

\newcommand\h{h}

\title{Explicit models for perverse sheaves}
\author{F. Gudiel Rod´r\'{\i}guez\thanks{Partially supported by BFM2001-3207 and TMR-40
``Singularities of Differential Equations and Foliations".} and L.
Narv\'aez Macarro\thanks{Supported by BFM2001-3207.}\\
Departamento de \'Algebra, Universidad de Sevilla}
\date{}
\begin{document}
\maketitle

\begin{abstract}
We consider categories of generalized perverse sheaves, with
relaxed constructibility conditions, by means of the process of
gluing $t$-struc\-tures and we exhibit explicit abelian categories
defined in terms of standard sheaves categories which are
equivalent to the former ones. In particular, we are able to
realize perverse sheaves categories as non full abelian
subcategories of the usual bounded complexes of sheaves
categories. Our methods use induction on perversities. In this
paper, we restrict ourselves to the two-strata case, but our
results extend to the general case.
\medskip

{\footnotesize \noindent{\sc Keywords:} perverse sheaf, derived
category, $t$-structure,
stratified space, abelian category. }
\end{abstract}

\section*{Introduction}

Perverse sheaves first appear in context of Complex Analytic
Geometry by the coming together of the Riemann-Hilbert
correspondence of Mebkhout-Kashiwara and the Intersection
Cohomology of Goresky-MacPherson at the beginning of the 1980s. In
the work \cite{bbd_83} the notion of $t$-structure over a
triangulated category was extracted  and it was proved that the
category of analytic constructible perverse sheaves, that we call
``classical perverse sheaves", can be obtained by a general
process of ``gluing" $t$-structures, that makes sense in a much
more general framework. In fact, the main contribution of {\it
loc.~cit.} is the use of that process to define  $\ell$-adic
perverse sheaves over algebraic varieties in positive
characteristics and to prove the theorem of  purity of
intersection complexes.
\medskip

In this paper, we develop some ideas  and complete some results in
 \cite{nar_lille} and \cite{gudiel_tesis} on the core of the
$t$-structure obtained by gluing standard $t$-structures shifted
by ``perversities" of strata, as in the classical case but without
imposing necessarily any constructibility conditions. Objects in
this core can be thought of as ``generalized perverse sheaves".

In the complex analytic case, and when we consider the middle
perversity, the category of classical perverse sheaves is a full
abelian subcategory of that of generalized perverse sheaves.
Furthermore, a classical perverse sheaf is the same as a
generalized perverse sheaf which is complex analytic
constructible.

The advantage of our point of view consists of being able to work
simultaneously with different perversities and to establish some
precise relations between perverse sheaves with respect to
different perversities, which we do not know how to do if we are
restricted to the classical case.

Our main result is theorem \ref{teo}, from which we deduce (see
\ref{apli-explicit}) that any (generalized) perverse sheaf, and
then any classical perverse sheaf, has a canonical model
(\ref{model}). As a consequence, the category of (generalized)
perverse sheaves is equivalent to a non full (resp. full) abelian
subcategory of the category of the usual bounded complexes (resp.
up to homotopy).

The main idea consists of constructing a functor $\Phi$ relating
$d$-perverse and $(d-1)$-perverse sheaves, and by iteration,
$d$-perverse sheaves with $0$-perverse sheaves, which are nothing
but usual sheaves. In this way we develop the idea pointed out in
\cite{nar_lille}, rem. (2.3.7), where we were restricted to the
``conical" case.

Construction of functor $\Phi$ and many other results in this
paper are inspired by the formalism of vanishing cycles
\cite{sga_7_II} and the gluing of classical perverse sheaves of
Deligne-Verdier \cite{del_81,ver_85} and MacPherson-Vilonen
\cite{mac_vi_86}, but our framework is more general.

In order to simplify, in this paper we restrict ourselves to the
two-strata case, but our results extend to the general case.

Let us now comment on the content of this paper.
\medskip

In section 1 we recall first the gluing process of $t$-structures
and the notion of (generalized) perverse sheaf is introduced.
Second, we recall some elementary constructions with adjoint
functors that play a fundamental role in the proof of theorem
\ref{teo} and in the manipulation of our explicit models for
perverse sheaves.
\medskip

Section 2 deals with the construction of functor $\Phi$ and the
``induction on perversities" \ref{induct}.
\medskip

In section 3 we show the main result in this paper, namely that
the category of (generalized) perverse sheaves $\Perv$ is
equivalent to an explicit abelian category described in terms of
abelian categories of usual sheaves.
\medskip

In section 4 we give some applications of theorem \ref{teo}.
First, we associate to any perverse sheaf a canonical model. More
precisely, we lift the inclusion functor of $\Perv$ into the
derived category to a faithful exact functor into the category of
usual bounded complexes.

Second, we lift the inclusion functor of $\Perv$ into the derived
category to a fully faithful functor into the category of bounded
complexes up to homotopy, $\komp^b$. In particular, $\Perv$ can be
realized as a full abelian subcategory of $\komp^b$.

Third, we give quiver descriptions of conical perverse sheaves
with respect to a $K(\pi,1)$ basis.

Finally, we compute in terms of our listed canonical models the
different perverse direct images and the intersection complex
associated to a sheaf on the open stratum, and we announce further
results.
\medskip

Part of the work of the first author was carried out during a
visit to CMAF da Universidade de Lisboa, whose hospitality is
gratefully acknowledged.

\section{Preliminaries and notations}

\subsection{Perverse sheaves} \label{preli}

Let $X$ be a topological space stratified by $\Sigma = \{C,U\}$,
where $i: C \rightarrow X$ is a closed immersion and
$j:U=X-C\rightarrow X$ is its complementary dense open immersion.
Let $\OO_X$
 be a sheaf of rings on $X$ and let $\OO_U = j^* \OO_X, \OO_C =
 i^*\OO_X$. For $*= X, U, C$, let us denote by $\gB_*$ the
 abelian category of sheaves of $\OO_*$-modules, and let
 $\gA_*\subset \gB_*$ be full abelian subcategories stable for
 kernels, cokernels and
extensions. Let us denote by $\calD_* :=\cd{+}{\gA_*}{\gB_*}$ the
full triangulated subcategory of the derived category
$\catder(\gB_*)$ whose objects are bounded below complexes with
cohomology in $\gA_*$. Let us suppose that the usual functors
$i_*=i_! , i^*, \RR i^!, \RR j_*, j_!, j^* = j^! $ induce functors
$$ \calD_C
\genfrac{}{}{0pt}{1}{\xrightarrow{i_*=i_!}}{\xleftarrow[i^*, \RR
i^!]{}}
 \calD_X
\genfrac{}{}{0pt}{1}{\xrightarrow{j^*=j^!}}{\xleftarrow[\RR j_*,
j_!]{}} \calD_U $$ in such a way that we are in the conditions of
gluing $t$-structures on $\calD_U$ and on $\calD_C$ \cite{bbd_83}.

\begin{ejemplo} \label{ejem-1}\\
(1)  If $\gA_* = \gB_* $, then $\calD_* = \cd{+}{}{\gB_* } $.\\
(2) Let $S$ be a compact topological space, $X$ the cone of $S$,
$C$ its vertex, $\OO_X$ the constant sheaf with fiber a ring
(resp. a noetherian ring) $k$ and $\gA_*$ the abelian categories
of $\Sigma$-constructible sheaves of $k$-modules not necessarily
finitely generated (resp. finitely generated).\\
(3) The space $X$ is a pseudomanifold stratified by $\Sigma$,
$\OO_X$ is the constant sheaf with a field $k$ as fiber, and the
$\gA_*$ are the abelian categories of $\Sigma$-constructible
sheaves of $k$-vector spaces of arbitrary (resp. of finite) rank.
For instance, $X$ can be a complex analytic space and $C \subset X$
a smooth closed analytic set satisfying the Whitney conditions.
\end{ejemplo}

\begin{definicion}{\label{perv}} For any integer $d \geq 0$,
the category of {\em $d$-perverse sheaves} on $X$ with respect to
the stratification $\Sigma $, $\Perv^d (X, \Sigma )$, is the core
of the $t$-structure on $\calD_X$ obtained by gluing the natural
$t$-structure on $\calD_U$ and the image by $[-d]$ of the natural
$t$-structure on $\calD_C$ \cite{bbd_83}. We will say that the
{\em perversity} of the stratum $C$ (resp. $U$) is $d$ (resp.
$0$).
\end{definicion}

\begin{nota}
Observe that, if $d=0$, then the category $\Perv^0 (X,\Sigma)$
coincides with the category $\gA_X$.
\end{nota}

\begin{proposicion} \label{carac}
{\em (Characterization of $d$-perverse sheaves)} An
object $K$ of  $\calD_X$ is a $d$-perverse sheaf (with respect to
$\Sigma$) if and only if the following properties hold:
\begin{enumerate}
\item[(a)]
$K$ is concentrated in degrees $[0, d]$, \item[(b)] $j^* K $ is
concentrated in degree $0$, \item[(c)] $h^n \RR i^! K = 0 $
 for $n<d$.
 \end{enumerate}
\end{proposicion}

\begin{demostracion}
By definition of $\Perv^d (X , \Sigma )$, an object $K$ of
$\calD_X$ is a $d$-perverse sheaf if and only if $h^n (j^* K) = 0
$ for $n \neq 0$, $h^n i^* K = 0 $ for $n>d$ and $h^n \RR i^! K = 0 $
for $n<d$, and it is clear that a $K$ satisfying  properties (a),
(b), (c) is $d$-perverse.

Let us now take  a $d$-perverse sheaf $K$. Properties (b) and (c)
are clear. The long exact sequence associated with the triangle
$$\tri{j_! j^* K}{K}{i_* i^* K}$$gives rise to isomorphisms $h^l
(K) \simeq h^l (i_* i^* K)$ for any $l\geq 1$ and then $h^l(K)=0$
for any $l>d$.

In a similar way, the long exact sequence associated with the
triangle $$\tri{i_* \RR i^! K}{K}{\RR j_* j^* K}$$ and the fact that $\RR
j_* j^* K $ is concentrated in non-negative degrees gives rise to
isomorphisms $h^l (i_* \RR i^! K)\simeq h^l (K) $ for $l<0$ and then
$h^l(K)=0$ for any $l<0$, and $K$ is concentrated in degrees
$[0,d]$.
\end{demostracion}

\subsection{Functors acting on morphisms of functors}

Let $\cB,\cC$ be categories, $F,G:\cB\to\cC$ functors and
$\tau:F\to G$ a morphism of functors (or natural transformation)
 which associates to any object $B$ in $\cB$ a morphism
$\tau_B:FB \longrightarrow GB$ in $\cC$ with the usual naturality
properties.
\medskip

For any functors $E:\cA\to\cB, H:\cC\to\calD$ we denote by $\tau
E: FE\to GE$, $H\tau: HF\to HG$ the morphisms given by $$ (\tau
E)_A = \tau_{EA},\quad (H\tau)_B = H(\tau_B)$$for any objects $A$
in $\cA$ and $B$ in $\cB$.
\medskip

\numero \label{rules} We have the following rules:
\begin{enumerate}\item[(a)] $ H(\tau E) = (H\tau)E$, $H 1_F
=1_{HF}$, $1_F E=1_{FE}$.
\item[(b)]
$(\tau\circ\varepsilon)E = (\tau E)\circ(\varepsilon E)$, $
H(\tau\circ\varepsilon) = (H\tau)\circ (H\varepsilon)$ for any
other morphism $\varepsilon:F'\to F$.
\item[(c)] $(\sigma G)\circ (K\tau) = (L\tau)\circ
(\sigma F)$ for any other functors $K,L:\cC\to\calD$ and any other
morphism $\sigma:K\to L$.\item[(d)] $ (\tau +\tau')E = (\tau
E)+(\tau' E),\quad H(\tau+\tau') = (H\tau)+ (H\tau')$ in the case
of additive functors between additive categories.
\end{enumerate}

\subsection{Adjoint functors}
\label{adjoint}

\numero {\label{adj}} In this section,  we consider a couple of
adjoint additive functors $G: \cA\to \cA'$, $F : \cA' \to \cA$
between abelian categories with adjunction morphisms $\alpha :
\Id_\cA \to FG $ and $\beta : GF \to \Id_{\cA'} $  such that $F$
is left exact, $G$ is exact and $\alpha $ is injective. We  denote
$ \F := F  G $ and $(\Q,q):= \coker \alpha$. We have then a
commutative diagram with exact rows and columns:
\begin{equation}\label{CD-1}
\begin{CD} {} @. 0 @.0 @. 0\\
@. @VVV @VVV @VVV\\
 0 @>>> \Id_\cA @>{\alpha}>> \F @>{q}>> \Q @>>> 0\\ @.
@V{\alpha}VV @V{\alpha\F}VV @V{\alpha\Q}VV @.\\ 0 @>>> \F
@>{\F\alpha}>> \F\F @>{\F q}>> \F\Q \\ @. @V{q}VV @V{q\F}VV
@V{q\Q}VV\\ 0 @>>> \Q @>{\Q\alpha}>> \Q\F @>{\Q q}>> \Q\Q \\ @.
@VVV @VVV @VVV\\{} @. 0 @.0 @. 0
\end{CD}\end{equation}
\medskip

\numero\label{gamma} Let us call $\gamma : \QQ \lra \F \F $ the
unique morphism satisfying $\gamma \circ q = \alpha\F - \F \alpha
$. From (\ref{CD-1}) we deduce the relations $$(\F q) \circ \gamma
= \alpha\QQ,\quad (q\F)\circ\gamma = - \Q\alpha.$$

\numero \label{mu} From the adjunction properties, the exact
sequence $$\suex{G}{G\alpha}{G\F}{Gq}{G\Q}$$splits, with
retraction $\beta G : GFG \lra G$. Then the sequence
\begin{equation}\label{ex-seq-1}\suex{\F}{\F\alpha}{\F\F}{\F q}{\F\Q}
\end{equation}is exact and splits, with retraction $\nu:= F\beta G$. Let us
call $\mu:\F\Q\to \F\F$ the corresponding section, i.e. $$
\mu\circ (\F q) = 1_{\F^2} - (\F\alpha)\circ\nu,\quad (\F
q)\circ\mu = 1_{\F\Q}.
$$
\medskip

\numero \label{retr-2}With the above notations, the relation
 $\nu \circ (\alpha\F) = 1_\F $
holds and the sequence
$$\suex{\F}{\alpha\F}{\F\F}{q\F}{\Q\F}$$also splits with the same
retraction as in (\ref{ex-seq-1}). Let  us call $\mu':\Q\F\to
\F\F$ the corresponding section, i.e. $$ \mu'\circ (q\F) =
1_{\F^2} - (\alpha\F)\circ\nu,\quad (q\F)\circ\mu' = 1_{\Q\F}. $$
We have $\gamma = \mu \circ (\alpha \Q)$.

Functors $\F\Q$ and $\Q\F$ are canonically isomorphic by means of
$\h:= (q\F)\circ\mu:\F\Q\to \Q\F$ and its inverse $\h^{-1}=(\F
q)\circ\mu'$.
\medskip

\begin{lema} \label{le:adj} For any objects $A,B$ in $\cal A$, the sequence
$$ \suex{\Hom (\Q A,\F B)}{q_A^*}{\Hom (\F A,\F
B)}{\alpha_A^*}{\Hom (A,\F B)}$$is exact and splits.
\end{lema}

\begin{demostracion} From \ref{retr-2}, application
$$f\in \Hom (A,\F B) \mapsto \nu_B \circ (\F f)\in \Hom (\F A,\F
B)$$ is a section of the above sequence.
\end{demostracion}

\section{Construction of categories and  functors}

\subsection{The functor $\Omega $}

\numero Let $\gA$ be a category. Let us denote by $\Arr(\gA)$ the
category of arrows of $\gA$, by $s, t: \Arr(\gA ) \rightarrow \gA$
the functors defined by $$s(A\xrightarrow{u} B):= A,\quad
t(A\xrightarrow{u} B):= B$$ and by $\zeta : s\to t$ the morphism
defined by $\zeta_{(A\xrightarrow{u} B)} := u$.
\medskip

If $\gA$ is abelian, the category $\Arr(\gA)$ is also abelian and
functors $s, t$ are exact and induce exact functors $\os,\ot:\comp
(\Arr (\gA )) \to \comp (\gA)$. They induce triangulated
 functors $ \komp (\Arr (\gA )) \to \komp
(\gA)$, $\catder(\Arr(\gA)) \to \catder(\gA)$, also denoted by
$\os,\ot$. Let us denote by $\ozeta:\os\to \ot$
 the morphism of functors induced by
$\zeta$.
\medskip

\numero \label{N}The functor $N:\comp(\Arr(\gA)) \to
\Arr(\comp(\gA))$ defined by $N=\os\xrightarrow{\ozeta}\ot$ is an
isomorphism of abelian categories. In a similar way we define
functors $N:\komp(\Arr(\gA)) \to \Arr(\komp(\gA))$,
$N:\catder(\Arr(\gA)) \to \Arr(\catder(\gA))$, which are no longer
equivalence of categories. Nevertheless, a morphism in $\komp
(\Arr (\gA ))$ is a quasi-isomorphism if and only if its images by
$\os$ and $\ot$ are quasi-isomorphisms.
\medskip

For any abelian category $\gA$ and any object
$(U\xrightarrow{\beta}V) \in \comp ( \Arr (\gA))\equiv \Arr(\comp
(\gA))$ we define $$\Omega (U \xrightarrow{\beta} V ) := ( V
\xrightarrow{q} \cone (\beta ))\in \comp ( \Arr (\gA)),$$where $q$
is the canonical inclusion. One can easily define the action of
$\Omega$ on morphisms and we obtain an exact functor $\Omega:
\comp ( \Arr (\gA))\to \comp ( \Arr (\gA))$ which commutes (up to
isomorphism) with the translation functor and satisfies $\os
\Omega = \ot$.

\begin{proposicion} The functor $\Omega$ above induces a
triangulated functor $\Omk: \komp ( \Arr (\gA))\to \komp ( \Arr
(\gA))$ such that $\os \Omk = \ot$.
\end{proposicion}

\begin{demostracion} It is an exercise we leave to the reader.
\end{demostracion}

The definition of distinguished triangles in $\komp(\gA)$ gives
rise to a morphism $\vartheta:\ot\Omk \to \os[1]$ in such a way
that the following triangle of functors
\begin{equation}\label{tri-0-k}
\trito{\os}{\ozeta}{\ot=\os \Omk}{\ozeta \Omk}{\ot
\Omk}{\vartheta}
\end{equation}
is distinguished, i.e. its evaluation on any object of $\komp (
\Arr (\gA))$ is a distinguished triangle of $\komp(\gA)$.
\medskip

The following proposition is basically the same as the axiom (TR2)
of triangulated categories for $\komp(\gA)$ (\cite{ver_cd},
chap.~I, prop.~3.3.3).

\begin{proposicion} \label{prop-133}Under the above hypothesis, there is an
isomorphism of functors $\chi: \os[1] \xrightarrow{\simeq} \ot
\Omk^2$ such that the following diagram is commutative:
$$\begin{CD}
\ot \Omk @>{\vartheta}>> \os[1]
@>{-\ozeta[1]}>> \ot[1]\\
@V{=}VV @V{\chi}V{\simeq}V @V{=}VV\\
\os\Omk^2 @>{\ozeta\Omk^2}>> \ot \Omk^2 @>{\vartheta\Omk}>>
\os[1]\Omk.
\end{CD}
$$
\end{proposicion}

\begin{proposicion} \label{prop-134} The functor $\Omk: \komp ( \Arr (\gA))\to \komp (
\Arr (\gA))$ transforms quasi-isomorphisms into quasi-isomorphisms
and then it induces a triangulated functor $\Omd: \catder ( \Arr
(\gA))\to \catder ( \Arr (\gA))$ such that $\os \Omd = \ot$.
Moreover, there is a morphism $\vartheta:\ot\Omd \to \os[1]$ and
an isomorphism $\chi: \os[1] \xrightarrow{\simeq} \ot \Omd^2$ such
that the following triangle of functors
\begin{equation}\label{tri-0-d}
\trito{\os}{\ozeta}{\ot=\os \Omd}{\ozeta \Omd}{\ot
\Omd}{\vartheta}
\end{equation}
is distinguished and the following diagram is commutative:
$$\begin{CD}
\ot \Omd @>{\vartheta}>> \os[1]
@>{-\ozeta[1]}>> \ot[1]\\
@V{=}VV @V{\chi}V{\simeq}V @V{=}VV\\
\os\Omd^2 @>{\ozeta\Omd^2}>> \ot \Omd^2 @>{\vartheta\Omd}>>
\os[1]\Omd.
\end{CD}
$$
\end{proposicion}

\begin{demostracion} The first part follows from the relation $\os\Omk =
\ot$, from the fact that a morphism $\xi $ in $\komp ( \Arr
(\gA))$ is a quasi-isomorphism if and only if $\os(\xi),\ot(\xi)$
are quasi-isomorphisms and from triangle (\ref{tri-0-k}).

The second part is basically the axiom (TR2) for the triangulated
category $\catder(\gA)$ and follows from triangle (\ref{tri-0-k})
and from proposition \ref{prop-133}.
\end{demostracion}

\begin{nota} The functor $\Omega_{\gA}$ defined in
\cite{nar_88,nar_lille} is related to the functor $\Omd$ above by
the equality $\Omega_{\gA} = N \Omd$.
\end{nota}

\numero \label{frakQ} Let us denote by $\frak Q$ the full
subcategory of $\Arr(\catder(\gA))$ whose objects are the
$A\xrightarrow{\upsilon}B$ such that $A$ is concentrated in degree
$0$ and $B$ is concentrated in degrees $\geq 0$. For such objects
the morphism $\upsilon$ (in $\catder(\gA)$) is determined by its
cohomology of degree $0$. More precisely we have the following
result:

\begin{proposicion} \label{cono}
Functor $N$ defines an equivalence of (additive) categories
between $N^{-1}{\frak Q}$ and $\frak Q$.
\end{proposicion}

\begin{demostracion} We sketch the definition of a quasi-inverse
of $N:N^{-1}{\frak Q}\to{\frak Q}$ and leave the details to the
reader.

Given an object $Y = (A\xrightarrow{\upsilon}B)$ in $\frak Q$, let
$U,V$ be the complexes defined by $U^0 = h^0A$, $U^n=0$ for $n\neq
0$ and $$V = \tau_{\geq 0}B = \cdots \xrightarrow{} 0
\xrightarrow{} \coker d_B^{-1} \xrightarrow{\overline{d_B^0}} B^1
\xrightarrow{d^1_B} \cdots,$$ where $\coker d_B^{-1}$ is placed in
degree 0, and let $\widetilde{\upsilon}:U\to V$ the morphism of
complexes determined by $\widetilde{\upsilon}^0 =
h^0{\upsilon}:U^0 \to h^0 B \subset V^0$.
\medskip

Correspondence $Y\mapsto (U\xrightarrow{\widetilde{\upsilon}} V)$
extends to a functor $\overline{N}:{\frak Q}\to N^{-1}{\frak Q}$.
It is easy to see that $N \overline{N}\simeq \Id_{\frak Q}$.
\medskip

On the other hand, for any object $X = (U \xrightarrow{\beta} V)$
in $N^{-1}{\frak Q}$, the following commutative diagram in $\comp
(\gA)$
$$\begin{CD}
U @>{\beta}>> V\\
@V{\text{nat.}}VV @VV{\text{nat.}}V \\
\tau_{\geq 0}U @>{\tau_{\geq 0}\beta}>> \tau_{\geq 0}V\\
@A{\text{nat.}}AA @AA{=}A\\
h^0 U @>{\widetilde{\beta}}>> \tau_{\geq 0} V\end{CD}$$ defines a
natural isomorphism
\begin{equation}\label{kappa}
\kappa(X): X \longrightarrow \overline{N}(N(X)).
\end{equation}
\end{demostracion}

\numero \label{propo} Let us call $C= (\ot\Omd
\overline{N})[-1]:{\frak Q}\to \catder(\gA)$.

An object $(U\xrightarrow{\beta}V)\in \catder(\Arr (\gA))$ is in
$\frak P:= \Omd^{-1}N^{-1}{\frak Q}$ if and only if the complex
$U$ is concentrated in degrees $\geq 0$, $V$ is concentrated in
degree $0$ and $h^0\beta$ is injective.

From propositions \ref{prop-134} and \ref{cono} we obtain an
isomorphism
\begin{equation}\label{eq-N}
\eta := (\ot\Omd \kappa \Omd)[-1]\circ \chi[-1]:
\os\longrightarrow C N\Omd
\end{equation}
between functors from $\frak P$ to $\catder (\gA)$.
\medskip

\numero \label{precision} For any object $Y =
(A\xrightarrow{\upsilon}B)\in \frak Q$ such that $A^n= 0 $ for all
$n\neq 0$ and $B^n =0$ for all $n < 0$ we can identify (by a
canonical isomorphism) \begin{equation} \label{preci-0} C(Y) =
\cdots\longrightarrow 0 \longrightarrow A^0
\xrightarrow{-\upsilon^0} B^0 \xrightarrow{-d^0_B} B^1
\xrightarrow{-d^1_B}\cdots,\end{equation}where $A^0$ is placed in
degree $0$.

Consequently, for any object $X=(U \xrightarrow{\beta} V)\in
{\frak P}$ such that $U$ and $V$ are concentrated in degree $0$
($h^0\beta$ must be injective), we can identify $$ (C N\Omd)(X) =
\cdots \xrightarrow{} 0 \xrightarrow{} h^0V
\xrightarrow{-\text{nat.}} \coker h^0\beta \xrightarrow{} 0
\xrightarrow{}\cdots$$ placed in degrees $0,1$ and isomorphism
$\eta_X:U\to (C N\Omd)(X)$ reduces to
\begin{equation}\label{preci-1}\begin{CD} h^0U @>{-h^0\beta}>> h^0V\\ @V{0}VV
@VV{-\text{nat.}}V\\ h^1U=0 @>{0}>> \coker h^0\beta,
\end{CD}\end{equation}
where minus sign in $-h^0\beta$ comes from the definition of
$\chi$ in proposition \ref{prop-134}.

\subsection{The functor $\Phi$}
\label{Phi}

 In this section we come back to the situation
described in section \ref{preli}.
\medskip

\numero \label{nume} Let us choose an additive left exact functor
$\F =FG : \cA={\gB}_U \rightarrow \cA={\gB}_U $ and an injective
morphism $\alpha:1\to \F$ as in \ref{adj}, such that $\F
(\gA_U)\subset \gA_U$, $ \F_{|\gA_U}$ is exact,
 $(\RR^i j_*) (\F A) = 0$ and
$j_* \F A \simeq  \RR (j_* \F)  A $, for $i>0, A \in {\gA}_U$. The
restriction to $\gA_U$ of the functor $\Q= \coker \alpha$ defined
in \ref{adj} is also exact.

To simplify, let us write $\Omega: \catder(\Arr (\gB_X)) \to
\catder(\Arr (\gB_X))$ instead of $\Omd$ in proposition
\ref{prop-134}.
\medskip

\numero \label{Psi}Let us first consider the additive left exact
functor $\psi_{\F} : {\gB}_X \rightarrow \Arr({\gB}_X) $ defined
by
\begin{equation}\label{rho}
\psi_{\F} = (\Id \xrightarrow{\rho} j_* \F j^*  ),\quad \rho :=
(j_* \alpha j^*)\circ \adj ,\end{equation}where $\adj : \Id\to
j_*j^*$ is the adjunction morphism, and second, functors
$$\Psi_{\F} := \Omega  \RR \psi_{\F} : \calD_X \rightarrow
\catder(\Arr (\gB_X)),\quad \Phi_{\F} := \ot \Psi_{\F}: \calD_X
\to \calD_X.$$ Once the functor $\F$ is fixed, we omit subscripts
and we will write $\psi,\Psi,\Phi$ instead of
$\psi_\F,\Psi_\F,\Phi_\F$.

 We have canonical isomorphisms
 \begin{equation}\label{cano-1}
\overline{s} \RR \psi  \simeq \Id,\quad  \os \Psi \simeq
\overline{t} \RR \psi  \simeq \RR (j_* \F j^*)= \RR (j_* \F)
j^*\end{equation} and triangle (\ref{tri-0-d}) gives rise to the
triangle
\begin{equation} \label{trian-2}
\trit{\Id}{\rho}{\RR (j_* \F) j^*}{u^1}{\Phi }.
\end{equation}

\begin{ejemplo}\label{ejem-2}\\
(1) In example \ref{ejem-1}, (1) let us take $U^{\text{dis}}$ as
the discrete topological space with underlying set $U$, $\Delta:
U^{\text{dis}} \to U$ the identity map, $\OO_{U^{\text{dis}}} =
\Delta^* \OO_U$, $\cA'$ the abelian category of
$\OO_{U^{\text{dis}}}$-modules and $F = \Delta_*, G = \Delta^*$
(see \cite{god_58}, chap.~II, \S 4.3).
\\
(2) In example \ref{ejem-1}, (2), let us suppose that $S$ is a
$K(\pi,1)$ space, $p:\widetilde{U}\to U$ the universal covering
space of $U$, $\OO_{\widetilde{U}} = p^* \OO_U$, $\cA'$ the
abelian category of $\OO_{\widetilde{U}}$-modules and $F = p_*, G
= p^*$ (see \cite{nar_lille}). If the fundamental group $\pi_1 (U,
x_0) $ is finite and $k$ is noetherian, then we can also consider
the categories $\gA_* $ as those  of constructible sheaves of
finitely generated modules.\\ (3) In example \ref{ejem-1}, (3), it
is not possible in general to choose a functor $\F$ as above, but
we will be able to apply the methods of this paper as explained in
section \ref{apli-explicit}.
\end{ejemplo}

\subsection{Induction on perversities}

 The following theorem generalizes \cite{nar_lille}, prop. 2.3.3 and rem. 2.3.7.

\begin{teorema} \label{induct} Let $d$ be an integer $\geq 1$ and let
 $ K$ be an object in $\calD_X$.
Then, $K \in \Perv^d(X,\Sigma) $ if and only if $ j^* K \in \gA_U$
and $\Phi K \in \Perv^{d-1}(X,\Sigma)$.\end{teorema}

\begin{demostracion} Let us consider the long exact sequence of cohomology
 associated with the
triangle (\ref{trian-2}) evaluated on $K$:
\begin{equation} \label{trian-2-1}
\trit{K}{\rho_K}{\RR (j_* \F )j^*K}{u^1_K}{\Phi K}.
\end{equation}

If $K$ is $d$-perverse, $j^*K$ belongs to $\gA_U$, $\RR(j_* \F)
j^* K = j_*\F j^* K$ is concentrated in degree $0$ and $ h^i (\Phi
K ) \simeq h^{i+1} (K ) $ for $ i \neq 0,-1$. In particular $\Phi
K$ is concentrated in degrees $[-1,d-1]$.

For $i=-1$ we have an exact sequence $$ \suext{h^{-1} (\Phi  K
)}{h^0  K }{\RR^0(j_* \F) j^* K= j_*\F j^* (h^0 K)},$$ where the
second arrow is nothing but
 $\rho_{h^0 K}$ (see (\ref{rho})),
which is injective because $\alpha$ is injective and
$\RR^0\Gamma_C(h^0 K) = \RR^0\Gamma_C K = i_* h^0 \RR i^! K =0$.
Then $h^{-1} (\Phi  K)=0$ and $\Phi  K$ is concentrated in degrees
$[0,d-1]$.

By applying $j^*$ to (\ref{trian-2-1}) we obtain $j^*\Phi  K
\simeq \Q j^* K$, and then $j^*\Phi  K$ is concentrated in degree
$0$.

On the other hand, by applying $\RR i^!$ to (\ref{trian-2-1}) we
deduce that $\RR i^!\Phi  K\simeq (\RR i^!K)[1]$, hence $h^m \RR i^!\Phi K =0$
for any $m<d-1$.

By proposition \ref{carac}, we conclude that $\Phi  K$ is
$(d-1)$-perverse.
\medskip

Conversely, let us suppose that $j^*K\in \gA_U$ and that $\Phi  K$
is $(d-1)$-perverse. By triangle (\ref{trian-2-1}) again we deduce
that $K$ is concentrated in degrees $[0,d]$, $h^m \RR i^! K \simeq
h^{m-1} \RR i^!\Phi  K =0$ for any $m<d$ and then $K$ is $d$-perverse.
\end{demostracion}

\begin{nota} As pointed out in \cite{nar_lille}, rem. (2.3.7), theorem \ref{induct}
suggests iterating the functor $ \Phi $ in order to obtain, for
any $d$-perverse sheaf $K$, an usual sheaf $\Phi^d K$. The main
result of this paper (see theorem \ref{teo}) tells us how to
reconstruct $K$ from its restriction to the open set $U$ and from
$\Phi^d K$.
\end{nota}

\section{The equivalence of categories}

\subsection{Gluing data}\label{catefun}

We keep the notations in \ref{adj} and \ref{nume}.
\medskip

\numero \label{g}For each integer $i\geq 0$ let us write $q^i := q
{\Q^i}, \alpha^i := \alpha{\Q^i}, g^i := (\F q^i)\circ
(\alpha\F\Q^i) = [(\F q)\circ (\alpha \F)]\Q^i = [(\alpha\Q)\circ
q]\Q^i = \alpha^{i+1}\circ q^i = g^0 \Q^i$.

For  $d\geq 1$, we have the complex
\begin{equation}\label{trunca}T_d := \F
\xrightarrow{g^0 } \cdots \xrightarrow{g^{d-2}} \F \Q^{d-1}
\xrightarrow{q^{d-1}}\Q^{d}\end{equation}placed in degrees
$[0,d]$, which is a resolution of length $d$ of the identity
functor by means of the injection $\Id \xrightarrow{\alpha}\F$.
\medskip

\numero \label{hd} Let $\h_d: T_{d}\Q\to \Q T_{d+1}$ be the
morphism of complexes given by $$ h_d^i = (-1)^{i+1}\h\Q^i:
\F\Q^{i+1}\to \Q\F\Q^i,\quad 0\leq i\leq d-1$$and $h_d^{d}=
(-1)^{d}\Q\alpha\Q^{d}:\Q^{d+1}\to \Q\F\Q^{d}$, where $h$ has been
defined in \ref{retr-2}. It is a quasi-isomorphism whose
composition with $\alpha\Q:\Q\to T_{d}\Q$ is equal to $\Q\alpha:
\Q\to \Q T_{d+1}$.
\medskip

\begin{definicion} \label{def-Bd}For each integer $d\geq 1$,
let $\gB_{\F}^d$ be the additive category whose objects are the
$({\calL} , {\calF}, u, \sigma )$, where ${\calL} \in \gA_U,\calF
\in \gA_X$, $ u: j_* \F \Q^{d-1}\calL \rightarrow\calF $  and
$\sigma : \Q^d\calL \xrightarrow{\sim} j^*\calF$ such that
 $u \circ j_* g^{d-2}_{\calL} = 0 $ (if $d\geq 2$) and $\sigma\circ q^{d-1}_{\calL} = j^*u$,
 and whose morphisms from $({\calL} , {\calF}, u, \sigma )$ to
 $({\calL}' , {\calF}', u', \sigma' )$  are defined as the pairs
 $(f,g)$ where $f:\calL\to \calL', g:\calF\to\calF'$ are morphisms
 such that $u'\circ (j_*\F\Q^{d-1} f) = g\circ u$.
\end{definicion}

Observe that in the above definition, the relation $(j^*g)\circ
\sigma = \sigma'\circ (\Q^d f)$ holds.
\medskip

For $d\geq 2$, let us define the additive functor $$\G_{d-1} :=
\coker j_* g^{d-2} : \gA_U \rightarrow \gA_X ,$$ which is right
exact.
\medskip

In the above definition, we can replace condition $u \circ j_*
g^{d-2}_{\calL} = 0 $ by taking objects $({\calL} , {\calF},
\overline{u}, \sigma )$ with $\overline{u}: \G_{d-1}\calL \to
\calF$.

The proof of the following proposition is an exercise left up to
the reader.

\begin{proposicion} The category $\gB_{\F}^d$ is abelian.
\end{proposicion}

\begin{nota}
By using the fact that sheaves on $X$ are determined by their
restrictions to $U$ and $C$ and by the gluing morphism $i^* \to
i^*j_*j^*$, category $\gB_{\F}^d$ fits into the construction of
abelian categories in \cite{mac_vi_86}, 1. Namely, category
$\gB_{\F}^d$ is equivalent to the category ${\scr C}(F, G; T)$ in
{\it loc.~cit.}, where $F= i^* \G_{d-1}, G= i^* j_* \Q^d : \gB_U
\rightarrow \gB_C $ ($F$ is right exact and $G$ is left exact) and
$T = i^* \overline{j_* q^{d-1}}:F\to G$, where $\overline{j_*
q^{d-1}}:\G_{d-1}\to j_*\Q^d$ is the morphism induced by $j_*
q^{d-1}$.

As in  \cite{mac_vi_86}, any other choice of the functor $\F$ in
\ref{nume} gives rise to a category equivalent to $\gB_{\F}^d$.
\end{nota}
\bigskip

\numero {\label{induc}} By theorem \ref{induct}, functors $j^*$
and $\Phi $ can be considered as functors
$j^*:\Perv^d(X,\Sigma)\to \gA_U,\quad \Phi : \Perv^d(X,\Sigma)\to
\Perv^{d-1}(X,\Sigma)$.
\medskip

From the properties of $\F$ in \ref{nume}, we have $$(j_*\F
)(\gA_U)\subset \bigcap_{m\geq 0} \Perv^m(X,\Sigma).$$

For any $K\in \Perv^d(X,\Sigma)$, we have $\RR (j_*\F )j^*K=j_*\F
j^*K$ and the morphism $u^1$ in (\ref{trian-2}) gives rise to a
morphism $$u^1: j_*\F j^* \to \Phi $$between functors $j_*\F j^*,
\Phi: \Perv^d(X,\Sigma)\to \Perv^{d-1}(X,\Sigma)$.
\medskip

As pointed out in the proof of theorem \ref{induct}, by applying
the functor $j^*$ to (\ref{trian-2}) we deduce an isomorphism
\begin{equation}\label{rel-1}
\xi^1: \Q j^* \xrightarrow{\sim} j^*\Phi \quad \text{such
that}\quad \xi^1\circ (q{j^*}) = j^*u^1.\end{equation}
\medskip

\numero \label{defi-0}We define inductively
\begin{eqnarray*}
u^{i}&=(u^{i-1}{\Phi }) \circ (j_*\F\Q^{i-2}\xi^1): j_*\F\Q^{i-1}
j^* \to \Phi ^i,\quad i\geq 2\\  \xi^i &= (\xi^1{\Phi ^{i-1}})
\circ (\Q \xi^{i-1}): \Q^i j^* \xrightarrow{\sim} j^* \Phi^{i},
\quad i\geq 1.\end{eqnarray*} The relations
\begin{equation}\label{rel-2}\xi^i \circ (q\Q^{i-1}j^*) = j^* u^i,\quad
u^i\circ (j_* g^{i-2}{j^*}) = 0,\quad (\xi^{i-1}\Phi)\circ
(\Q^{i-1}\xi^1) = \xi^i
\end{equation} hold for every $i\geq 2$.

\subsection{The theorem}

With the notations introduced in \ref{defi-0}, we do the
following:
\medskip

\begin{definicion} \label{defi-Bd}For each integer $d\geq 1$,
let us define the additive functors $$D_d := (j^*, \Phi ^d, u^d,
\xi^d) : \Perv^d(X,\Sigma) \to \gB_{\F}^d$$  and $B_d : \gB_{\F}^d
\rightarrow \Perv^d (X,\Sigma)$ by
\begin{equation} \label{model}
B_d (\calL ,\calF , u , \sigma ):= j_* \F\calL \xrightarrow{j_*
g^0_\calL} j_* \F \Q\calL \xrightarrow{j_* g^1_\calL } \cdots
\xrightarrow{j_* g^{d-2}_\calL } j_* \F \Q^{d-1}\calL
\xrightarrow{u}\calF,
\end{equation}
where the complex is placed in degrees $[0,d]$, the action of
$B_d$ on morphisms being obvious.
\end{definicion}

In the above definition we can identify
\begin{equation}\label{identi}
j^*B_d (\calL,\calF , u , \sigma ) = T_d\calL\end{equation}by
means of $\sigma$ (see (\ref{trunca})). Furthermore, the
acyclicity properties in \ref{nume} show that $j_*\F\Q^i\calL =
\RR j_*\F\Q^i\calL$. Then $\RR i^!B_d (\calL ,\calF , u , \sigma )
= \RR i^!\calF[-d]$, and  we deduce the perversity of $B_d (\calL
,\calF , u , \sigma )$ from proposition \ref{carac}.
\medskip

The main result of this paper is the following:
\medskip

\begin{teorema}\label{teo} For any integer $d\geq 1$,
functors $B_d$ and $D_d$ defined above are the quasi-inverse of
each other and they define, thus, an equivalence of categories
between $\Perv^d (X,\Sigma) $ and $\gB_{\F}^d$.
\end{teorema}

As suggested by rem. (2.3.7) of \cite{nar_lille} and theorem
\ref{induct}, the proof of theorem \ref{teo} can be approached  by
induction on perversity $d$.

\begin{nota}
In case $d=1$ our proof of the isomorphism $\Id \simeq D_1 B_1$ is
essentially the same as in \cite{nar_lille}, th. 2.3.4, but it
should be noticed that in {\it loc.~cit.} there is a mistake in
the proof of the faithfulness of $D_1$. Our proof of theorem
\ref{teo} completes the one given in \cite{gudiel_tesis}.
\end{nota}
\medskip

\subsection{The proof}

\noindent {\bf First Part:}  We are going to construct a natural
isomorphism $\calO \simeq D_d B_d\calO$ for any $\calO$ in
$\gB_{\F}^d$. \bigskip

 For $d=0$ let us call $\gB_{\F}^0 = \gA_X$ and $B_0
=\Id : \gB_{\F}^0 \to \Perv^0(X,\Sigma)$.
\medskip

For any $d\geq 2$ we consider functors $\gT : \gB^d_\F \to
\gB^{d-1}_\F $ whose action on objects (resp. on morphisms) is
given by $ \gT (\calL , \calF , u , \sigma ) := (\Q \calL , \calF
, u , \sigma )$ (resp. $\gT(f,g) := (\Q f,g)$). For $d=1$, functor
$\gT: \gB^1_\F \to \gB^{0}_\F $ is simply defined by $\gT (\calL ,
\calF , u , \sigma ) = \calF$.
\medskip

For any $d\geq 1$ we also consider functors $s: \gB^d_\F \to
\gA_U$, $t: \gB^d_\F \to \gA_X$ and morphism $\upsilon:
j_*\F\Q^{d-1} s \to t$ defined by $s(\calL , \calF , u , \sigma )
= \calL$, $t(\calL , \calF , u , \sigma )=\calF$,
$\upsilon_{(\calL , \calF , u , \sigma )} = u$. We obviously have
$s\gT = \Q s$ and $j^* t \stackrel{\sigma}{=} \Q^d s$.
\medskip

From (\ref{identi}) we can identify $j^* B_d = T_d s$ for $d\geq
1$, and from the acyclicity properties of $\F$ with respect to
$j_*$ in \ref{nume}, we deduce
$$\Phi B_d = \ot \Omega \RR
\psi B_d  = \ot \Omega \psi B_d = \cone(B_d \xrightarrow{\rho B_d
} j_*\F j^* B_d = j_*\F T_d s),$$i.e.
$$
(\Phi B_d)^{-1} = 0\oplus (j_* \F s),\quad (\Phi B_d)^{d-1} = (j_*
\F^2\Q^{d-1} s)\oplus t,$$
$$(\Phi B_d)^{d} = (j_*
\F\Q^{d}s)\oplus 0,$$ $$(\Phi B_d)^{i} = (j_* \F^2\Q^i s)\oplus
(j_*\F\Q^{i+1}s), \quad 0\leq i\leq d-2,$$ and
$$ d^{-1}_{\Phi B_d} =
\left(\begin{array}{cc} 0&j_*\alpha{\F s}\\
0& -j_* g^0 s\end{array}\right),\quad d^{i}_{\Phi B_d} =
\left(\begin{array}{cc} j_* \F g^i s & j_*\alpha {\F\Q^{i+1} s}\\
0& -j_*g^{i+1} s\end{array}\right), \quad 0\leq i\leq d-3,$$
$$d^{d-2}_{\Phi B_d} = \left(\begin{array}{cc} j_* \F g^{d-2} s &
j_*\alpha{\F\Q^{d-1}s}\\ 0& -\upsilon
\end{array}\right),\quad
d^{d-1}_{\Phi B_d} = \left(\begin{array}{cc} j_* \F q {\Q^{d-1} s}
& (j_*\alpha {\Q^{d} s})\circ \adj\\ 0& 0\end{array}\right),$$
where $\adj:t\to j_* j^* t = j_* \Q^d s$ is the adjunction
morphism.
\medskip

By the same reason, morphism $u^1 B_d:j_*\F j^*B_d=j_*\F T_d s \to
\Phi B_d$ becomes the natural inclusion.
\medskip

Let $Q_d$ be the complex of functors from $\gB^d_{\F}$ to $\gA_U$
obtained by plumbing $\F s$ in degree $-1$ and $\F T_d s$ in
degrees $\geq 0$ by means of $\F\alpha s$. From (\ref{trunca}) and
\ref{nume} we deduce that complexes $Q_d$ and $j_*Q_d$ are exact.
\bigskip

\numero \label{lambdas} For any $d\geq 1$ let us call $L_d =
B_{d-1}\gT$, which can be considered as a complex of functors from
$\gB^d_{\F}$ to $\gA_X$, and let us define the following exact
sequences: $$\textstyle\suex{L_d^{-1}=0}{\lambda_d^{-1}=0}{(\Phi
B_d)^{-1}}{\pi_d^{-1}= (0\ 1)}{j_*Q_d^{-1}},$$
$$\textstyle\suex{L_d^i}{\lambda_d^i
=(-1)^{i+1}\left(\scriptsize\begin{array}{c} j_* \mu {\Q^i s}\\
-1\end{array}\right)}{(\Phi B_d)^i}{\pi_d^i= ( 1 \ j_*\mu {\Q^i
s})}{j_*Q_d^i},\quad 0\leq i\leq d-2,$$
$$\textstyle\suex{L_d^{d-1}} {\lambda_d^{d-1} = (-1)^{d}\left(
\substack{(j_*\gamma {\Q^{d-1} s })\circ (\adj) \\-1} \right)}
{(\Phi B_d)^{d-1}} {\pi_d^{d-1} = ( 1 \ (j_*\gamma {\Q^{d-1}
s})\circ (\adj) )} {j_*Q_d^{d-1}},$$
$$\textstyle\suex{L_d^d=0}{\lambda_d^d=0}{(\Phi
B_d)^d}{\pi_d^d=(1\  0)}{j_* Q_d^d},$$ where $\gamma:\Q\to \F^2$,
$\mu:\F\Q\to \F^2$ have been defined in \ref{gamma} and \ref{mu}
respectively.
\medskip

From \ref{gamma}, \ref{mu}, \ref{g} and \ref{retr-2} we deduce,
first: $$(\F g^i)\circ(\mu\Q^i) = [(\F\alpha\Q)\circ(\F
q)\circ\mu]\Q^i = \F\alpha\Q^{i+1},$$ $$(\mu\Q^{i})\circ g^{i} =
[\mu\circ(\F q)\circ (\alpha\F)]\Q^{i} = [\alpha \F
-\F\alpha]\Q^{i},$$ $$(j_*\gamma {\Q^{d-1} s})\circ (\adj)\circ
\upsilon =(j_*\gamma {\Q^{d-1} s})\circ (j_* q {\Q^{d-1} s}) =
[j_* (\alpha\F-\F\alpha)\Q^{d-1}]s,$$ $$(\F q {\Q^{d-1}s})\circ
(\gamma \Q^{d-1}s) - \alpha \Q^{d} s= [(\F q)\circ \gamma -
\alpha\Q] {\Q^{d-1} s} = 0,$$ and second: $$d^i_{\Phi
B_d}\circ\lambda_d^i = \lambda_d^{i+1}\circ (j_*g^{i+1} s),\quad
d^i_{Q_d}\circ \pi_d^i = \pi_d^{i+1}\circ d^{i+1}_{\Phi
B_d}\quad\text{for any $i$.} $$ In particular we obtain an exact
sequence of complexes $$ \suex{B_{d-1}\gT}{\lambda_d}{\Phi
B_d}{\pi_d}{j_*Q_d},$$which shows that $\lambda_d:B_{d-1}\gT\to
\Phi B_d$ is a quasi-isomorphism and then an isomorphism between
functors from $\gB^d_\F$ to $\Perv^{d-1}(X,\Sigma)$.
\bigskip

\numero  \label{ci-1}For any $d\geq 1$ we consider the morphism of
functors $\theta_d:j_*\F s \to B_{d-1}\gT$ given by $\theta_d =
j_*g^0 s$ if $d\geq 2$ and $\theta_1 =\upsilon: j_*\F s\to t =
B_0\gT$. Diagram
\begin{equation} \label{cd-1}
\begin{CD} j_*\F s
@>{\theta_d}>> B_{d-1}\gT \\ @V{j_*\F \alpha s}VV
@V{\lambda_d}VV\\ j_*\F j^*B_d = j_*\F T_d s @>{u^1 B_d}>> \Phi
B_d\end{CD}\end{equation}commutes in the homotopy category of
complexes and then in the derived category.
\bigskip

\numero \label{ci-2} For each $i\geq 1$ let us call $$\phi_i:=
(j_*\F\Q^{i-1}\xi^{-1} B_{i+1})\circ (j_*\F\Q^{i-1} j^*
\lambda_{i+1}) : j_*\F\Q^{i-1}j^* B_i \gT \longrightarrow
j_*\F\Q^i j^* B_{i+1}.$$

From rule \ref{rules}, (c),  and \ref{defi-0} we deduce that the
following diagram  of functors from $\gB^{i+1}_\F$ to
$\Perv^0(X,\Sigma)=\gA_X$ $$\begin{CD} j_*\F\Q^{i-1}j^* B_i \gT
@>{u^i B_i\gT}>> \Phi^i B_i \gT\\ @V{\phi_i}VV
@V{\Phi^i\lambda_{i+1}}VV\\ j_*\F\Q^i j^* B_{i+1}
@>{u^{i+1}B_{i+1}}>> \Phi^{i+1} B_{i+1}\end{CD}$$commutes, where
the vertical arrows are isomorphisms.
\medskip

\numero With identifications $$ j_*\F\Q^{i-1}j^* B_i \gT =
j_*\F\Q^{i-1}T_i s\gT = j_*\F\Q^{i-1}T_i \Q s$$and$$ j_*\F\Q^i j^*
B_{i+1} = j_*\F\Q^{i-1}\Q T_{i+1} s$$one can prove that $\phi_i =
j_*\F\Q^{i-1} h_i s$, where $h_i$ has been defined in \ref{hd},
but we will not need that result in the rest of this paper.
\bigskip

Summing up \ref{ci-1} and \ref{ci-2}, for any $d\geq 1$ we obtain
a commutative diagram of functors from $\gB^d_\F$ to
$\Perv^0(X,\Sigma)=\gA_X$ $$ \begin{CD} j_*\F s\gT^{d-1}
@>{\upsilon \gT^{d-1}}>> B_0 \gT\gT^{d-1} = e \gT^{d-1}\\
@V{j_*\F\alpha s \gT^{d-1}}VV  @VV{\lambda_1 \gT^{d-1}}V\\ j_*\F
j^* B_1 \gT^{d-1}  @>{u^1B_1 \gT^{d-1}}>> \Phi B_1 \gT^{d-1}\\
@V{\phi_1 \gT^{d-2}}VV @VV{\Phi\lambda_2\gT^{d-2}}V\\ \vdots  @.
\vdots\\ @V{\phi_{d-2} \gT}VV @VV{\Phi\lambda_{d-1}\gT}V\\
j_*\F\Q^{d-2} j^* B_{d-1} \gT @>{u^{d-1}B_{d-1} \gT}>> \Phi^{d-1}
B_{d-1}\gT\\ @V{\phi_{d-1} }VV @VV{\Phi\lambda_{d}}V\\
j_*\F\Q^{d-1} j^* B_{d} @>{u^{d}B_{d}}>> \Phi^{d}
B_{d}.\end{CD}$$Compositions of vertical arrows give rise to the
natural isomorphism $\Id_{\gB^d_\F} \xrightarrow{\sim} D_d B_d$ we
wanted and the first part of the proof of theorem \ref{teo} is
finished.
\bigskip

\noindent {\bf Second part:}

In this part we prove that for any $d$-perverse sheaf $K$, there
exists a natural isomorphism $K\simeq B_d D_d K$. We are using
notations of \ref{propo}. We proceed by induction on $d\geq 1$.
\medskip

For any $d$-perverse sheaf we know (theorem \ref{induct}) that
$\RR\psi K\in{\frak P}$ and $$ N \Psi K = j_*\F j^* K
\xrightarrow{u^1_K} \Phi K.$$

Let us call $$\omega_K: K\xrightarrow{\sim} C(j_*\F j^* K
\xrightarrow{u^1_K} \Phi K)$$the composition of isomorphism
$$\eta_{\RR\psi K}: \os \RR\psi K \to CN\Omega \RR\psi K = CN\Psi
K = C(j_*\F j^* K \xrightarrow{u^1_K} \Phi K)$$ defined in
\ref{propo} and isomorphism $K\xrightarrow{\sim} \os \RR\psi K$ of
(\ref{cano-1}).
\medskip

Functors $\psi,\Omega, N$ commute with $j^*$ and we can identify
\begin{equation} \label{*}
j^* C(j_*\F j^* K \xrightarrow{u^1_K} \Phi K) \stackrel{\xi_1}{=}
 C(\F j^* K \xrightarrow{q_{j^*K}} \Q j^* K).
\end{equation}

Then, by using (\ref{preci-1}) we obtain
\begin{equation}\label{**}
j^*\omega_K = -\alpha_{j^*K}: j^*K \longrightarrow C(\F j^* K
\xrightarrow{q_{j^*K}} \Q j^* K).\end{equation} \medskip

For $d=1$ we have $B_1 D_1 K = j_*\F j^* K \xrightarrow{u^1_K}
\Phi K$ which is isomorphic to $$C(j_*\F j^* K \xrightarrow{u^1_K}
\Phi K)= j_*\F j^* K \xrightarrow{-u^1_K} \Phi K$$ by means of
$(1, -1)$. The composition of this last isomorphism with
$\omega_K$ gives rise to an isomorphism $$ \delta^1_K: K
\longrightarrow B_1 D_1 K$$ natural with respect to $K\in
\Perv^1(X,\Sigma)$ such that $j^* \delta^1_K = -\alpha_{j^*K}$.
\medskip

Now let $d$ be an integer $\geq 2$ and suppose there exists
$\delta^{d-1}: \Id_{\Perv^{d-1}(X,\Sigma)}\xrightarrow{\sim}
B_{d-1} D_{d-1}$ such that
\begin{equation}\label{***} j^* \delta^{d-1} = (-1)^{d-1} \alpha
j^*: j^* \longrightarrow j^* B_{d-1} D_{d-1}= T_{d-1} s D_{d-1} =
T_{d-1} j^* .\end{equation} Isomorphism $$(\xi^1,1): \gT D_d = (\Q
j^*, \Phi^d, u^d, \xi^d) \to D_{d-1}\Phi = (j^* \Phi, \Phi^{d},
u^{d-1}\Phi, \xi^{d-1}\Phi)$$allows us to identify both functors
and, by (\ref{***}) we obtain $$ (j^*\delta^{d-1}\Phi)\circ (j^*
u^1) = (-1)^{d-1} g^0 j^*: \F j^* \longrightarrow j^* B_{d-1}
D_{d-1} \Phi = j^* B_{d-1} \gT D_d = T_{d-1} \Q j^*.$$Then
$(\delta^{d-1}\Phi)\circ u^1 = (-1)^{d-1} j_*g^0 j^*$ and $$
C(1,\delta^{d-1}\Phi): C(j_*\F j^* \xrightarrow{u^1} \Phi )
\xrightarrow{\sim} C(j_*\F j^* \xrightarrow{(-1)^{d-1}j_* g^0 j^*}
B_{d-1} \gT D_d),$$but
\begin{eqnarray*}
& B_{d-1} \gT D_d  =  j_*\F\Q j^*  \xrightarrow{j_*g^0{\Q j^* }}
\cdots \xrightarrow{j_*g^{d-3}{\Q j^* }} j_*\F\Q^{d-2}\Q
\xrightarrow{u^d} \Phi^{d-1}\Phi&
\end{eqnarray*}
and $j_*g^{i-1}{\Q j^* } = j_* g^{i}{j^*}$. In particular, by
using (\ref{preci-0})
 we deduce an isomorphism
\begin{eqnarray}
\nonumber &C(j_*\F j^* \xrightarrow{u^1} \Phi ) \simeq&\\
\label{c-0}& j_*\F \xrightarrow{(-1)^d j_*g^0j^*}  j_*\F\Q j^*
\xrightarrow{-j_*g^1j^* } \cdots \xrightarrow{-j_*g^{d-2}j^* }
j_*\F\Q^{d-1} \xrightarrow{-u^d} \Phi^{d},&
\end{eqnarray}
and the complex (\ref{c-0}) is isomorphic to $$B_d D_d= j_*\F
\xrightarrow{j_*g^0j^*}  j_*\F\Q j^* \xrightarrow{j_*g^1j^* }
\cdots \xrightarrow{j_*g^{d-2}j^* } j_*\F\Q^{d-1}
\xrightarrow{u^d} \Phi^{d}$$ by means of $$
((-1)^{d-1},-1,1,-1,\dots,(-1)^{d-1},(-1)^d).$$By composing
isomorphisms above with $\omega$ we obtain an isomorphism
$$\delta^{d}: \Id_{\Perv^{d}(X,\Sigma)}\xrightarrow{\sim} B_{d}
D_{d}$$ such that $$ j^* \delta^{d} = (-1)^{d-1}(- \alpha j^*)=
(-1)^d\alpha j^*$$and the proof of theorem \ref{teo} is finished.

\section{Applications}

\subsection{Explicit models for perverse sheaves}
\label{apli-explicit}

Theorem \ref{teo} provides explicit models (\ref{model}) for
$d$-perverse sheaves. Actually, functor $B_d$ factorizes through
the category of bounded complexes $\comp^b (\gB_X)$ and it defines
a faithful exact functor $B_d:\gB^d_\F \to \comp^b (\gB_X)$
establishing an equivalence of categories between $\gB^d_\F$ and a
non full abelian subcategory of $\comp^b (\gB_X)$, whose objects
are precisely complexes of the form (\ref{model}). In particular,
inclusion functor $\Perv^d(X,\Sigma)\subset \calD_X$ can be lifted
to an exact faithful functor $\Perv^d(X,\Sigma)\to\comp^b(\gB_X)$.
\medskip

The lifting above allows us to describe in a concrete way the
realization functor (see \cite{bbd_83}, 3.1.9)
 $$ \text{real} :
\catder(\Perv^d(X,\Sigma))\longrightarrow \calD_X$$by taking
single complexes associated with double complexes.
\medskip

When no functor $\F$ is available for the given subcategories
$\gA_*\subset \gB_*$, we can always work at the level of the full
derived categories $\catder^+(\gB_*)$ by using, for instance,
Godement functor $\F = \Delta_* \Delta^*$, as shown in examples
\ref{ejem-1}, (1) and \ref{ejem-2}, (1). The corresponding
category of perverse sheaves $\Perv^d(X,\Sigma)$ (without any
constructibility conditions, i.e. $\gA_* = \gB_*$) is, by theorem
\ref{teo}, equivalent to $\gB^d_\F$, whose objects are
\ref{def-Bd} the $(\calL,\calF,u,\sigma)$ where ${\calL} \in
\gB_U,\calF \in \gB_X$, $ u: j_* \F \Q^{d-1}\calL \rightarrow\calF
$  and $\sigma : \Q^d\calL \xrightarrow{\sim} j^*\calF$ such that
 $u \circ j_* g^{d-2}_{\calL} = 0 $ (if $d\geq 2$) and $\sigma\circ q^{d-1}_{\calL} = j^*u$,
 and whose morphisms from $({\calL} , {\calF}, u, \sigma )$ to
 $({\calL}' , {\calF}', u', \sigma' )$  are defined as the pairs
 $(f,g)$ where $f:\calL\to \calL', g:\calF\to\calF'$ are morphisms
 such that $u'\circ (j_*\F\Q^{d-1} f) = g\circ u$.
\medskip

Let us call $\Perv^d_c(X,\Sigma)$ the category of perverse sheaves
``constructible" with respect to $\gA_*\subset \gB_*$. It is a
full (abelian) subcategory of $\Perv^d(X,\Sigma)$ and then it is
equivalent to the full subcategory $\gB^d_{\F,c}$ of $\gB^d_\F$
whose objects are the $(\calL,\calF,u,\sigma)$ such that
$\calL\in\gA_U$ and morphism $\overline{u}:\G_{d-1}\calL \to
\calF$ has kernel and cokernel in $\gA_X$.
\medskip

So, even when no functor $\F$ is available for the given
subcategories $\gA_*\subset \gB_*$, explicit models and liftings
as above also exist.
\medskip

\begin{ejemplo} {\it (Perverse sheaves categories which split)}
In example \ref{ejem-1}, (2), let $d\geq 2$ be an integer and let
us suppose  $S$ a ``good"   compact, connected and simply
connected topological space, and $k$ a field such that
\begin{equation}\label{condi}\cohom^i(S,k)=0\quad \forall i=1,\dots,d.\end{equation}

For example, $S$ can be the
 $(n-1)$-dimensional sphere and $X$ the  $n$-dimensional disk,
stratified by the origin and its complement, for $n\geq d$, or in
singularity theory, $(X,0)\subset \C^{d+1}$ is an isolated
hypersurface singularity with complex link $S$ a topological
(exotic) sphere \cite{mil}, \S 8.
\medskip

Let us consider the category of $\Sigma$-constructible complexes
of sheaves of $k$-vector spaces    of arbitrary\footnote{We may
also consider only sheaves of finite rank.}
 rank on each
stratum which are $d$-perverse sheaves, denoted by
 $\Perv^d_c(X,\Sigma)$. It is a full  subcategory of
the category of $d$-perverse sheaves (without constructibility
conditions) $\Perv^d(X,\Sigma)$, which is equivalent by theorem
\ref{teo} to category $\gB^d_\F$, with $\F$ a functor satisfying
the  conditions  \ref{nume} (see \ref{ejem-2}, (1)).

Since  $S$  is simply connected, any locally-constant sheaf
$\calL$ of $k$-vector spaces on $U$ is constant with fiber $E=
\Gamma(U,\calL) \simeq k^r$ and $$(\RR^i j_* \calL)_C =
\lim_{\epsilon\to 0} \cohom^i(]0,\epsilon[\times S,\calL) =
\cohom^i(S,E)=0, \quad 1\leq i\leq d.$$ In particular, the
sequence $$ \suex{j_*\calL}{j_*\alpha_\calL}{j_*\F\calL}{j_*
q_\calL}{j_*\Q\calL}$$is exact and $\RR^i j_* \QQ\calL \simeq
\RR^{i+1} j_* \calL$ for all $i\geq 1$. Reasoning inductively we
obtain that the sequences
\begin{equation}\label{exac}
\suex{j_*\Q^{i-1}\calL}{j_*\alpha^{i-1}_\calL}{j_*\F\Q^{i-1}\calL}{j_*
q^{i-1}_\calL}{j_*\Q^i\calL},\quad i=1,\dots,d \end{equation}are
exact.
\medskip

Given a constructible  $d$-perverse sheaf
$K\in\Perv^d_c(X,\Sigma)$, let us denote $(\calL,\calF,u,\sigma) =
D_d K$ its corresponding object of $\gB^d_\F$ by means of theorem
\ref{teo}. Now $K$ is naturally isomorphic to
\begin{equation}\label{comple}j_* \F\calL \xrightarrow{j_*
g^0_\calL} j_* \F \Q\calL \xrightarrow{j_* g^1_\calL } \cdots
\xrightarrow{j_* g^{d-2}_\calL } j_* \F \Q^{d-1}\calL
\xrightarrow{u}\calF.\end{equation} The exactness of (\ref{exac})
for $i=d-1,d$ implies that $\coker j_* g^{d-2}_\calL =
j_*\Q^d\calL$. Let $s:j_*\Q^d\calL\to \calF$ be the morphism
induced by $u$, whose restriction to $U$ coincides with $\sigma$.
Now,  $\sigma$ being  an isomorphism, the
 adjunction properties for $(j^*,j_*)$ give us
a morphism $t:\calF\to j_*\Q^d\calL$ verifying  $ts= 1$. Then,
complex (\ref{comple}) is the  direct sum of $j_*T_d\calL$ and
 $(\ker t)[-d]$. On the other hand, the  exactness of (\ref{exac})
implies the $j_*T_d\calL$ is concentrated in degree $0$, its
$0$-cohomology  being equal to $j_*\calL$ and, thus, a
 constant sheaf. Finally, we obtain a natural isomorphism
$$ K\simeq (h^0K) \oplus (h^d K)[-d]$$expressing the
category $\Perv^d_c(X,\Sigma)$ as a direct sum of the
 category of constant sheaves of $k$-vector spaces in $X$
and the category of $k$-vector spaces, considered (this last
category) as the category of complexes of sheaves on $X$
 concentrated in degree $-d$ and supported by the vertex $C$.
\medskip

This is a  purely topological result related to a well-known
result of Kashiwa\-ra-Kawai \cite{kas_kaw_III} (see
\cite{mac_vi_86}, 6.5, p.~427). It can be also directly deduced by
using functors $j_!^p, j_*^p$ instead of our models.
Namely\footnote{We owe this remark to P. Deligne.}, our hypothesis
imply that $j_!^p j^* K \xrightarrow{\sim} j_*^p j^* K
\xrightarrow{\sim} j_* \calL \simeq h^0 K$ and then, from the
canonical morphisms
$$ j_!^p j^* K \xrightarrow{} K \xrightarrow{} j_*^p j^* K$$we
deduce that $h^0 K$ is a direct factor of $K$.
\end{ejemplo}

\subsection{Perverse sheaves categories as full abelian
subcategories of $\komp^b (\gB_X)$}

In this section we show that functor $B_d:\gB^d_\F \to
\komp^b_{\gA_X} (\gB_X) $ is fully faithful and then the inclusion
functor $\Perv^d(X,\Sigma)\subset\calD_X$ lifts to a fully
faithful functor $\Perv^d(X,\Sigma) \to \komp^b_{\gA_X} (\gB_X) $.
In particular, category $\Perv^d(X,\Sigma)$ is realized as a full
abelian subcategory of $\komp^b_{\gA_X} (\gB_X)$.
\medskip

\begin{teorema} \label{th:teo-2}Functor $B_d:\gB^d_\F \to
\komp^b_{\gA_X} (\gB_X)$ is fully faithful.
\end{teorema}

\begin{demostracion} Let $\calO_i = (\calL_i , \calF_i , u_i , \sigma_i
)$, $i=1,2$ be two objects in $\gB^d_\F$. We have to prove that $$
B_d:\Hom_{\gB^d_\F}(\calO_1,\calO_2) \xrightarrow{} \Hom_{\komp^b
(\gB_X)}(B_d\calO_1,B_d\calO_2)$$ is bijective.
\medskip

\noindent {\sc Injectivity:} Although the injectivity of $B_d$ is
a consequence of theorem \ref{teo} (the morphism $B_d:
\Hom_{\gB_{\F}^d} (\calO_1 , \calO_2 ) \to \Hom_{\calD_X} (B_d
(\calO_1 ) , B_d (\calO_2 )) $ is  bijective), we give here a
direct independent proof.

Let $(f,g):\calO_1\to\calO_2$ be a
morphism such that $B_d(f,g)$ is null-homotopic. We obviously have
$f= j^* h^0 B_d(f,g)=0$ and $B_d(0,g)$ is null-homotopic.
\medskip

There exist $s^i:j_*\F\Q^i\calL_1 \to j_*\F\Q^{i-1}\calL_2$,
$i=1,\dots,d-1$, $s^d:\calF_1\to j_*\F\Q^{d-1}\calL_2$ such that
$s^1\circ j_*g^0_{\calL_1} = 0$, $s^2\circ j_*g^1_{\calL_1} +
j_*g^0_{\calL_2}\circ s^1 =0$, $\dots$, $s^d\circ u_1 + j_*
g^{d-2}_{\calL_2}\circ s^{d-1}=0$, $g= u_2\circ s^d$,
$$\xymatrix{ j_* \F \calL_1 \ar[d]_0 \ar[r]^{j_* g^0_{\calL_1}}
 & j_* \F \Q \calL_1
\ar[dl]_{s^1}
 \ar[r]^{j_* g^1_{\calL_1}}  &  \dots &  j_*
\F \Q^{d-1} \calL_1 \ar[r]^{\qquad u_1} \ar[d]_0 \ar[dl]_{s^{d-1}}
&
 \calF_1 \ar[dl]_{s^d}
 \ar[d]_g \\
j_* \F \calL_2  \ar[r]_{j_* g^0_{\calL_2}} & \dots  & j_* \F
\Q^{d-2} \calL_2    \ar[r]_{j_* g^{d-2}_{\calL_1}} & j_* \F
\Q^{d-1} \calL_2 \ar[r]_{\qquad u_2}  & \calF_2 .}$$

In degree $0$, from $0=j^*s^1\circ g^0_{\calL_1} = j^*s^1\circ
\alpha^1_{\calL_1} \circ q_{\calL_1}$ we deduce $0=j^*s^1\circ
\alpha^1_{\calL_1}$ and then, there exists $t^1:\Q^2\calL_1\to
\F\calL_2$ s.t. $t^1\circ q^1_{\calL_1} = j^*s^1$. From  lemma
\ref{le:adj}, there exists $\tau^1: \F\Q^2\calL_1\to \F\calL_2$
s.t. $\tau^1\circ \alpha^2_{\calL_1}= t^1$.
\medskip

In degree $1$, from$$0 = j^* s^2\circ g^1_{\calL_1} +
g^0_{\calL_2}\circ j^* s^1 = (j^*s^2+ g^0_{\calL_2}\circ
\tau^1)\circ \alpha^2_{\calL_1} \circ q^1_{\calL_1}$$ we deduce
$$0=(j^*s^2+ g^0_{\calL_2}\circ \tau^1)\circ \alpha^2_{\calL_1}$$
and then, there exists $t^2:\Q^3\calL_1\to \F\Q\calL_2$ s.t.
$$j^*s^2+ g^0_{\calL_2}\circ \tau^1 = t^2 \circ
q^2_{\calL_1}.$$From  lemma \ref{le:adj} again, there exists
$\tau^2: \F\Q^3\calL_1\to \F\Q\calL_2$ s.t. $\tau^2\circ
\alpha^3_{\calL_1}= t^2$.
\medskip

We inductively construct $$t^i:\Q^{i+1}\calL_1\to
\F\Q^{i-1}\calL_2,\quad 2\leq i\leq d-1,$$
$$\tau^i:\F\Q^{i+1}\calL_1\to \F\Q^{i-1}\calL_2,\quad 2\leq i\leq
d-2$$ such that $$g^{i-2}_{\calL_2} \circ \tau^{i-1}+ j^*s^i = t^i
\circ q^i_{\calL_1},\quad 2\leq i\leq d-1,$$ $$ \tau^i\circ
\alpha^{i+1}_{\calL_1} = t^i,\quad 2\leq i\leq d-2. $$ Let us
identify $j^*\calF_1 = \Q^d\calL_1, j^*u_1 = q^{d-1}_{\calL_1}$ by
means of $\sigma_1$.
\medskip

In degree $d-1$, from $0=s^d\circ u_1 + j_* g^{d-2}_{\calL_2}\circ
s^{d-1}$ we deduce first
 \begin{multline*} 0= j^*s^d \circ q^{d-1}_{\calL_1} +
g^{d-2}_{\calL_2}\circ j^* s^{d-1} {=} j^*s^d \circ
q^{d-1}_{\calL_1} + g^{d-2}_{\calL_2}\circ (t^{d-1}\circ
q^{d-1}_{\calL_1} - g^{d-3}_{\calL_2}\circ \tau^{d-2}) = \\ =
(j^*s^d + g^{d-2}_{\calL_2}\circ t^{d-1})\circ
q^{d-1}_{\calL_1}\end{multline*}and second $$ 0 = j^*s^d +
g^{d-2}_{\calL_2}\circ t^{d-1}.$$ But $s^d$ is determined by its
restriction $j^*s^d$ $$ s^d = - (j_* g^{d-2}_{\calL_2})\circ (j_*
t^{d-1}) \circ (\text{adj}),$$where $\text{adj}: \calF_1 \to j_*
j^* \calF_1 = j_* \Q^d\calL_1$ is the adjunction morphism. Then$$
g = u_2\circ s^d = - u_2\circ (j_* g^{d-2}_{\calL_2})\circ (j_*
t^{d-1}) \circ (\text{adj}) = 0$$and injectivity is proven.
\bigskip

\noindent {\sc Surjectivity:} We need to prove that for any
morphism of complexes $F^\bullet: B_d \calO_1 \to B_d \calO_2$,
there exists $(f,g):\calO_1\to\calO_2$ s.t. $B_d(f,g)$ is
homotopic to $F^\bullet$.
\medskip

Obviously, morphism $f:\calL_1\to\calL_2$ must be equal to $j^*
h^0 F^\bullet$.
\medskip

Let us consider the following commutative diagram with exact
arrows$$
\begin{CD}
0 @>>> \calL_1 @>{\alpha_{\calL_1}}>> \F\calL_1 @>{
g^0_{\calL_1}}>> \F\Q\calL_1\\ @. @V{0}VV @V{j^*F^0-\F f}VV
@V{j^*F^1-\F\Q f}VV\\ 0 @>>> \calL_2 @>{\alpha_{\calL_2}}>>
\F\calL_2 @>{g^0_{\calL_2}}>> \F\Q\calL_2.\end{CD}$$ There exists
$\sigma^1_0:\Q\calL_1\to \F\calL_2$ s.t. $\sigma^1_0 \circ
q_{\calL_1} = j^*F^0-\F f$. From lemma \ref{le:adj}, there exists
$\sigma^1: \F\Q\calL_1\to \F\calL_2$ s.t. $\sigma^1 \circ
\alpha^1_{\calL_1} = \sigma^1_0$, and then $\sigma^1 \circ
g^0_{\calL_1} = j^*F^0-\F f$. Writing $s^1 := j_* \sigma^1$, we
have $s^1 \circ j_* g^0_{\calL_1} = F^0 - j_* \F f$.
\medskip

In a similar way, we inductively construct $s^i: j_*\F\Q^i\calL_1
\to j_* \F\Q^{i-1}\calL_2$, $i=2,\dots, d-1$, s.t. $s^i \circ
j_*g^{i-1}_{\calL_1} + j_* g^{i-2}_{\calL_2}\circ s^{i-1} =
F^{i-1}-j_*\F\Q^{i-1} f$. Let us write $\sigma^i = j^* s^i$.
\medskip

In degree $d-1$ we have
\begin{multline*}
(j^*F^{d-1} -\F\Q^{d-1} f - g^{d-2}_{\calL_2}\circ
\sigma^{d-1})\circ \alpha^{d-1}_{\calL_1}\circ
q^{d-2}_{\calL_1}=\\ =(j^*F^{d-1} -\F\Q^{d-1} f -
g^{d-2}_{\calL_2}\circ \sigma^{d-1})\circ g^{d-2}_{\calL_1}=\\
=g^{d-2}_{\calL_2}\circ (j^*F^{d-2} -\F\Q^{d-2} f - j^*F^{d-2}
+\F\Q^{d-2} f - g^{d-3}_{\calL_2}\circ \sigma^{d-2}) =0
\end{multline*}and then
$$0= (j^*F^{d-1} -\F\Q^{d-1} f - g^{d-2}_{\calL_2}\circ
\sigma^{d-1})\circ \alpha^{d-1}_{\calL_1}.$$ There exists
$\sigma^d_0: \Q^d\calL_1 \to \F\Q^{d-1}\calL_2$ s.t. $$\sigma^d_0
\circ q^{d-1}_{\calL_1} = j^*F^{d-1} -\F\Q^{d-1} f -
g^{d-2}_{\calL_2}\circ \sigma^{d-1}.$$ Since $j^*\calF_1 =
\Q^d\calL_1$, morphism $\sigma^d_0$ determines another morphism
$s^d: \calF_1 \to j_*\F\Q^{d-1}\calL_2$ s.t. $$s^d \circ u_1 + j_*
g^{d-2}_{\calL_2}\circ s^{d-1} = F^{d-1} - j_*\F\Q^{d-1} f.$$ To
finish, we take $$g:= F^d -u_2\circ s^d: \calF_1\to \calF_2.$$ An
straightforward computation shows that $(f,g):\calO_1\to \calO_2$
is a morphism in $\gB^d_\F$, and clearly the $s^i$, $i=1,\dots,
d$, give an homotopy between $F^\bullet$ and $B_d(f,g)$.
\end{demostracion}

\begin{corolario} The inclusion functor
$\Perv^d(X,\Sigma)\subset\calD_X = \catder^+_{\gA_X}(\gB_X)$ lifts
to a fully faithful functor $\Perv^d(X,\Sigma) \to \komp^b_{\gA_X}
(\gB_X) $. In particular, category $\Perv^d(X,\Sigma)$ is realized
as a full abelian subcategory of $\komp^b_{\gA_X}
(\gB_X)$.\end{corolario}

\begin{demostracion} It is a direct consequence of theorems
\ref{teo} and \ref{th:teo-2}.
\end{demostracion}

\subsection{Conical perverse sheaves with respect to a $K(\pi,1)$
basis}

In case of examples \ref{ejem-1}, (2) and \ref{ejem-2}, (2), we
suppose that $S$ is connected and its universal covering space is
contractible. Let us choose a base point $x_0\in S$ and let us
denote  $H=\pi_1(S,x_0)=\pi_1(U,x_0)$. Let $\gA_U$ be (resp.
$\gA_X$) the abelian category of locally constant sheaves of
$k$-modules (not necessarily finitely generated) on $U$ (resp. of
$\Sigma$-constructible sheaves of $k$-modules on $X$). We can take
$\F = p_* p^*$, where $p$ is the universal covering space of
$(U,x_0)$.
\medskip

Objects of category $\Perv^d(X,\Sigma)\subset
\catder^+_{\gA_X}(k_X)$ are called ``conical perverse sheaves" in
\cite{nar_lille}, def.~(2.1.1) and rem.~(2.3.7).
\medskip

\numero\label{coni-1} The standard equivalence of categories
between $\gA_U$ and $\Mod(k[H])$
 allows us to translate the exact sequence of functors of  $\gA_U$ $$
\suex{\Id}{\alpha}{\F}{q}{\Q}$$ in the following way. For each
$k[H]$-module $E$ we have:\\ 1) $\F E = E^H= \{ f : H \to E \}$,
where the action of $H$ is given by $ (hf) ( \sigma ) = f(\sigma
h),\quad f\in E^H, h, \sigma \in H.$\\ 2) Adjunction morphism
$\alpha_E:E \to \F E$ is given by $(\alpha_E e) (\sigma ) = \sigma
e$, $e\in E, \sigma\in H.$\\ 3) $\Q E = \{ \psi : H \to E\ |\
\psi(1)=0 \}$,
 where the action of $H$ is $$ (h \psi ) (\sigma) = \psi(\sigma h) -
 \sigma \psi(h),\qquad \psi \in \Q E,  \sigma , h \in H. $$
4) Morphism $q_E:\F E \to \Q E$ is given by $$ (q_E f)(\sigma) =
f(\sigma) - \sigma f(1),\quad f\in \F E = E^H, \sigma\in H.$$5)
The application  $c: e\in E \mapsto c(e)\in E^H$, where
$c(e)(\sigma) = e$  for any $\sigma\in H$, gives rise to a natural
identification $E = (\F E)^{inv}$.\\
6) For any $r\geq 1$ we have a natural identification$$\Q^r E = \{
\psi : H^r\to E\ |\ \psi (h_1,\dots,h_r) = 0\ \text{\ if\ }\
\exists j, h_j=1\}$$where the action of $H$ is given by
\begin{eqnarray*} &(h_{r+1} \psi)(h_1,\dots,h_r) = \displaystyle\sum_{i=1}^r
(-1)^{r-i}
\psi(h_1,\dots,h_{i-1},h_ih_{i+1},h_{i+2},\dots,h_{r+1}) +&\\
&+(-1)^r h_1\psi(h_2,\dots,h_{r+1}).\end{eqnarray*} 7) Morphisms
 $q^r_E: \F\Q^r E \to \Q^{r+1}E$, $g^r_E: \F\Q^r E \to
\F\Q^{r+1}E $ (see \ref{g}) are given by
\begin{eqnarray*} & (q^r_E f ) (h_1 , \dots, h_r, h_{r+1} ) =
f(h_{r+1}) (h_1 , \dots, h_r) - [h_{r+1} f(1)] (h_1 , \dots,
h_r),&\\ &(g^r_E f)(\sigma) = \sigma (q^r_E f),\quad  f\in \F \Q^r
E = (\Q^r E)^H,\ h_1,\dots, h_{r+1},\sigma\in H.&
\end{eqnarray*}8) By 5), morphism
$\varrho^r_E := (g^r_E)^{inv}: (\F\Q^r E)^{inv}=\Q^r E \to
(\F\Q^{r+1} E)^{inv}=\Q^{r+1}$ is $$ (\varrho^r_E
\psi)(h_1,\dots,h_{r+1}) = \psi ( h_1, \dots , h_r )- (h_{r+1}
\psi)( h_1, \dots , h_r )$$for $r\geq 1$, $\psi \in \Q^{r} E $,
$h_i\in H$. For $r=0$, morphism $\varrho^0_E := (g^0_E)^{inv}:
(\F\Q^0 E)^{inv}=E \to (\F\Q^{1} E)^{inv}=\Q^{r+1}$ is
$$ (\varrho^0_E e)(h_1) = e - h_1 e,\quad e\in E, h_1\in H.$$

\begin{nota} The complex $(\Q^r E,\varrho^r_E)_{r\geq 0}$ is the
usual complex of $E$-valued cochains obtained from the normalized
bar resolution \cite{maclane-63}, chap.~IV, \S 5.
\end{nota}

\numero\label{coni-2} Category $\gA_X$ is equivalent to the
category\\ -) whose  objects are triplets $(V,W,\varsigma )$ where
$V$ is a $k[H]$-module (representing  the restriction $j^*$ of a
constructible sheaf), $W$ is a $k$-module (representing the fiber
$i^*$ at $F$  of a constructible sheaf) and $\varsigma:W\to
V^{inv}$ is a $k$-linear morphism (representing the adjunction
morphism $i^* \to i^*j_*j^*$).\\ -) whose morphisms are defined in
the obvious way.\medskip

By \ref{coni-1}, \ref{coni-2} and the fact that sheaves on $X$ are
determined by their restrictions $j^*, i^*$ and the adjunction
morphism $i^*\to i^*j_*j^*$, we deduce that category $\gB^d_{\F}$
is equivalent to the category $\gC^d(k,H)$:\\ -) whose objects are
4-uples  $(E,M,u,v)$ where $E$ is a $k[H]$-module, $M$ is a
$k$-module and $u,v$ appear in a commutative diagram $$ \xymatrix{
\Q^{d-2} E \ar[r]^{\varrho^{d-2}_E} &\Q^{d-1} E
\ar[rr]^{\varrho^{d-1}_E} \ar[dr]_{u} & & (\Q^d E)^{inv}
\\ & & M \ar[ur]_v }
$$such that $u \circ \varrho^{d-2}_E=0$, if $d\geq 2$.\\ -) whose
morphisms are defined in the obvious way.
\medskip

By theorem \ref{teo} we conclude that the category of $d$-conical
perverse sheaves is equivalent to $\gC^d(k,H)$.
\medskip

In case $d=1$, by defining $v_\sigma(y) = -v(y)(\sigma)$,
$\sigma\in H, y\in M$, we obtain an equivalence between
$\gC^1(k,H)$ and the category of $k$-module diagrams $$ E
\genfrac{}{}{0pt}{1}{\xrightarrow{\quad
 u \quad }}{\xleftarrow[ \{v_\sigma\}_{\sigma\in H} ]{}} M
 $$such that\\
 (1)
 $ v_{\tau \sigma} = v_\tau \circ u \circ v_\sigma + v_\tau +
 v_\sigma$ for all $\sigma, \tau\in H$.\\
 (2) $ 1_E + v_\sigma\circ u$ is an automorphism of $E$ for any
 $\sigma\in H$.
 \medskip

Property (1) comes from the fact that  $v(y) \in (\Q
 E)^{inv}$ for every $y\in M$. In property (2), automorphism $1_E + v_\sigma\circ u$
 coincides with the action of $\sigma$ on $E$.
 \medskip

In this way we find a new proof of theorem (2.3.4) in
\cite{nar_lille}. This theorem is a natural generalization of the
first known case \cite{del_81} on explicit description of perverse
sheaves, namely $S= S^1, H = \Z$. (see also \cite{mai_bex},
\cite{nar_88}).

\subsection{Explicit description of perverse direct images and intersection complexes}

In this section we give models (\ref{model}) for $j_*^p \calL$,
$j_!^p \calL$ and $j_{!*}\calL$, where $\calL$ is an object of
$\gA_U$. The computations consist of interpreting the proof of
theorem 1.4.10 in \cite{bbd_83} in terms of our
$(\F,\Q)$-resolutions (\ref{trunca}).
\medskip

\numero For each $\calL\in \gA_U$ we have a natural isomorphism $$
D_d( j_*^p \calL) \simeq (\calL,j_*\Q^d\calL,j_*q^d_\calL,1).$$ In
particular
 the complex $$j_* \F\calL \xrightarrow{j_* g^0_\calL} j_* \F \Q\calL
\xrightarrow{j_* g^1_\calL } \cdots \xrightarrow{j_* g^{d-2}_\calL
} j_* \F \Q^{d-1}\calL \xrightarrow{j_*q^d_\calL}j_*\Q^d\calL $$
(in degrees $[0,d]$) is an explicit model for $j_*^p \calL$, which
coincides with $\tau_{\leq d} \RR j_* \calL$ \cite{bbd_83},
prop.~1.4.23.
 \bigskip

\numero For $d=1$ we have a natural isomorphism $$ D_1(
 j_!^p \calL) \simeq
(\calL,j_*\F\calL/j_!\calL,\text{can},1).$$ In particular the
complex $$ j_*\F\calL \xrightarrow{\text{can}}
j_*\F\calL/j_!\calL$$(in degrees $0,1$) is  an explicit model for
$j_!^p \calL$. It is quasi-isomorphic to $j_!\calL$ since
$j_!\calL$ is $1$-perverse.
\medskip

For $d\geq 2$ we have a natural isomorphism $$ D_d(
 j_!^p \calL) \simeq (\calL,\coker j_*
g^{d-2}_\calL ,\text{can},1).$$ In particular  the complex $$
j_*\F\calL \xrightarrow{j_*g^0_\calL}\cdots
\xrightarrow{j_*g^{d-2}_\calL} j_*\F\Q^{d-1}\calL
\xrightarrow{\text{can}} \coker j_* g^{d-2}_\calL,$$(in degrees
$[0,d]$) is an explicit model for $ j_!^p \calL$, which coincides
with $\tau_{\leq d-2} \RR j_* \calL$ \cite{bbd_83}, prop.~1.4.23.
\bigskip

\numero By interpreting natural morphisms $j_!^p \calL
\longrightarrow  j_*^p \calL$ on models above, we have a natural
isomorphism $$ D_d(j_{!*}\calL)\simeq (\calL,\Img j_*
q_{\Q^{d-1}\calL},j_* q_{\Q^{d-1}\calL},1).$$ In particular the
complex $$ j_*\F\calL \xrightarrow{j_*g^0_\calL}\cdots
\xrightarrow{j_*g^{d-2}_\calL} j_*\F\Q^{d-1}\calL \xrightarrow{j_*
q_{\Q^{d-1}\calL}} \Img j_* q_{\Q^{d-1}\calL} $$(in degrees
$[0,d]$) is an explicit model for the intersection complex
$IC(\calL) = j_{!*}\calL$, which coincides with $\tau_{\leq d-1}
\RR j_* \calL$ \cite{bbd_83}, prop.~1.4.23.

\subsection{Further results}

Following a suggestion of Deligne, explicit models of perverse
sheaves can be constructed by using other functorial resolutions
instead of (\ref{trunca}). For instance, given $\F =FG :
\cA={\gB}_U \rightarrow \cA={\gB}_U $, $\alpha:1\to \F$ under the
conditions of \ref{adj}, with $\F^k \calL$ $j_*$-acyclic for
$k\geq 1$ and $\calL \in {\gA}_U$, and not requiring $\F
({\gA}_U)\subset {\gA}_U$, we can use the ``simplicial" resolution
$$ \F \xrightarrow{ \partial^0} \F^2
\xrightarrow{\partial^1} \cdots \xrightarrow{
\partial^{d-2}} \F^d \xrightarrow{\partial^{d-1}} \cdots$$
where
$$ \partial^i = \alpha \F^{i+1} - \F \alpha \F^i + \cdots +
(-1)^{i+1} \F^{i+1} \alpha$$ (cf. \cite{god_58}, Appendice, 5 and
\cite{maclane-cwm}, VII, 6). This is the aim of an article in
preparation.

\bigskip

\noindent Departamento de \'Algebra\\
Facultad de Matem\'aticas, Univ. Sevilla\\
P.O. Box 1160\\
41080 Sevilla\\
SPAIN
\bigskip

\noindent E-mail:  gudiel@algebra.us.es, narvaez@algebra.us.es


\begin{thebibliography}{10}

\bibitem{bbd_83}
A.A. Beilinson, J.~Bernstein, and P.~Deligne.
\newblock {\em Faisceaux pervers}, volume 100 of {\em Ast\'erisque}.
\newblock S.M.F., Paris, 1983.

\bibitem{del_81}
P.~Deligne.
\newblock Lettre \`a {R.} {M}ac{P}herson.
\newblock 1981.

\bibitem{sga_7_II}
P.~Deligne and N.~Katz.
\newblock {\em Groupes de monodromie en G\'eom\'etrie Alg\'ebrique {(SGA 7
  II)}}, volume 340 of {\em Lect. Notes in Math.}
\newblock Springer-Verlag, Berlin-Heidelberg, 1973.

\bibitem{god_58}
R.~Godement.
\newblock {\em Topologie Alg\'ebrique et Th\'eorie des faisceaux}, volume XIII
  of {\em Publ. de l'Inst. de Math. de l'Univ. de Strasbourg}.
\newblock Hermann, Paris, 1958.

\bibitem{gudiel_tesis}
F.~Gudiel-{Rodr\'{\i}guez}.
\newblock Descripci{\'o}n expl{\'{\i}}cita de {$t$}-estructuras sobre espacios
  estratificados.
\newblock Univ. Sevilla, February 2001.
\newblock Ph.D.

\bibitem{kas_kaw_III} M.~Kashiwara and T.~Kawai.
\newblock On the holonomic systems of microdifferential equations {III}.
\newblock {\em Publ. Res. Inst. Math. Sci.} 17, (1981) 813--979.

\bibitem{maclane-63}
S.~Mac~Lane.
\newblock {\em Homology}.
\newblock Academic Press Inc., Publishers, New York, 1963.

\bibitem{maclane-cwm}
S.~Mac~Lane.
\newblock {\em Categories for the working mathematician}.
\newblock Springer-Verlag, New York, 1971.
\newblock Graduate Texts in Mathematics, Vol. 5.

\bibitem{mac_vi_86}
R.~MacPherson and K.~Vilonen.
\newblock Elementary construction of perverse sheaves.
\newblock {\em Invent. Math.} 84, (1986) 403--436.

\bibitem{mai_bex}
Ph. Maisonobe.
\newblock Faisceaux pervers sur {$\mathbb{C}$} relativement \`a {$\{0\}$} et couple
  {$E \substack{\rightarrow\\ \leftarrow} F$}.
\newblock In {L}\^e {D}\~ung {T}r\'ang, editor, {\em Introduction \`a la
  th\'eorie alg\'ebrique des syst\`emes diff\'erentiels}, pages 135--146.
  Hermann, Paris, 1988.
\newblock (Travaux en cours, vol. 34).

\bibitem{mil}
J.~Milnor.
\newblock {\em Singular points of complex hypersurfaces}, volume~61 of {\em
  Ann. of Math. Studies}.
\newblock Princeton Univ. Press, Princeton, N.J., 1968.

\bibitem{nar_88}
L.~Narv\'aez-Macarro.
\newblock Cycles \'evanescents et faisceaux pervers: cas des courbes planes
  irr\'eductibles.
\newblock {\em Comp. Math.} 65, (1988) 321--347.

\bibitem{nar_lille}
L.~Narv{\'a}ez-Macarro.
\newblock Cycles \'evanescents et faisceaux pervers. {I}{I}. {C}as des courbes
  planes r\'eductibles.
\newblock In J.P. Brasselet, editor, {\em Singularities (Lille, 1991)}, volume
  201 of {\em London Math. Soc. Lecture Note Ser.}, pages 285--323. Cambridge
  Univ. Press, Cambridge, 1994.

\bibitem{ver_85}
{J.L.} Verdier.
\newblock Extension of a perverse sheaf over a closed subspace.
\newblock {\em Ast\'erisque} 130, (1985) 210--217.

\bibitem{ver_cd}
J.L. Verdier.
\newblock Des cat\'egories d\'eriv\'ees des cat\'egories ab\'eliennes.
\newblock {\em Ast\'erisque} 239, 1996.
\newblock With a preface by Luc Illusie, Edited and with a note by Georges
  Maltsiniotis.
\end{thebibliography}
\end{document}